\documentclass[11pt,twoside]{article}

\usepackage{mathrsfs,amsfonts,amsmath,amssymb}
\usepackage{latexsym}
\usepackage{enumerate}
\usepackage[cp1250]{inputenc}
\usepackage{authblk}
\newtheorem{satz}{Theorem}
\newtheorem{proposition}[satz]{Proposition}
\newtheorem{theorem}[satz]{Theorem}
\newtheorem{lemma}[satz]{Lemma}

\makeatletter
\@namedef{subjclassname@2010}{%
  \textup{2010} Mathematics Subject Classification}
\makeatother

\def\no{\noindent}
\def\sbeq{\subseteq}

\def\Z{\mathbb {Z}}
\def\E{\mathbb {E}}

\def\e{\varepsilon}
\def\l{\lambda}
\def\s{\sigma}

\def\a{\alpha}

\def\C{{\mathbb C}}
\def\P{\mathbb {P}}
\def\h{\widehat}
\def\G{\Gamma}
\def\g{\gamma}
\def\d{\delta}

\def\({\big (}
\def\){\big )}
\def\g{\gamma}

\def\l{\lambda}

\def\m{\mu}

\def\n{\nu}

\def\r{{\rho}}

\def\zn{\Z/N\Z}

\def\v{\textrm{Var}}

\title{\bf On convex equations}
\author{{\sc Tomasz Schoen}}
\affil{Faculty of Mathematics and Computer Science,\\ Adam Mickiewicz
University,\\ Uniwersytetu Pozna\'nskiego 4, 61-614 Pozna\'n, Poland\\ schoen@amu.edu.pl}
\date{}
\begin{document}
\maketitle
\begin{abstract} 
We prove that every subset of $\{1,\dots, N\}$ which  does not contain any solutions to the equation $x+y+z=3w$ has
at most $$ {\exp\(-c(\log N)^{1/5+o(1)}}\)N$$ elements, for some $c>0$.  This theorem improves upon previous estimates. Additionally, our method has the potential to yield an optimal estimate for this problem that matches the known Behrend's lower estimate. Our approach relies on a new result on almost-periodicity of convolutions.
\end{abstract}



\bigskip
\centerline{\sc { 1. Introduction}}
\bigskip
The problem of estimating the maximum size of  a subset of $\{1,\dots, N\}$ that does not contain solutions to an invariant  linear equation 
$$a_1x_1+\dots+a_kx_k=0,$$
where $a_1, \dots, a_k\in \Z, k\ge 3$ and $a_1+\dots+a_k=0$, appears to be one of the most intriguing in additive combinatorics. In the general case not much is known, and even for the 'simplest' equation \cite{ruzsa}
$$2x_1+2x_2=y_1+3y_2$$
we are only able to prove only an upper bound of the form $N^{1-o(1)}$ and the lower bound $c N^{1/2}.$ More precise estimates are known only in the case of certain  symmetric equations i.e.
$$a_1x_1+\dots+a_kx_k=a_1y_{1}+\dots+a_ky_{k},$$
and equations of the form 
\begin{equation}\label{convex-eq}
a_1x_1+\dots+a_kx_k=by,
\end{equation}where $a_1,\dots, a_k, b$ are positive integers and $a_1+\dots+a_k=b,$ which we call convex.
We focus here on  the later type of equations. An important issue  is the fact that  we have a  general lower bound of a similar shape that applies to all convex equations. Behrend \cite{behrend} showed  (see also \cite{elkin} and \cite{green-wolf} for slight improvements) that for any convex equation \eqref{convex-eq} there is  a solution free subset of 
$\{1,\dots, N\}$ of size
$${\exp\(-C(\log N)^{1/2})\)}N,$$
where $C>0$ and it depends only on the coefficients of the given equation. The most interesting convex equation is, of course,
$x+y=2z$, as each of its nontrivial solutions (different form $x=y=z$) forms a nontrivial three term arithmetic progression. 
The first significant upper estimate was established by Roth \cite{roth} ( Roth's theorem, along with all the results quoted below, provides the same bounds for any invariant equation with an equal or greater number of variables).

\begin{theorem}(Roth)
Suppose that $A\sbeq \{1,\dots, N\}$ does not contain any nontrivial three term arithmetic progression. Then
$$|A|\ll \frac{N}{\log\log N}.$$
\end{theorem} 
Despite considerable effort and a high level of interest in the problem, we were only able to improve Roth's theorem by replacing his upper bound with  $N/(\log N)^{1+c},$
for some small constant $c>0$. Very recently, Kelley and Meka \cite{km} (see also \cite{bs}) proved a breakthrough result, much improving previous bounds, which is quite close to Behrend's lower bound.

\begin{theorem}(Kelley-Meka)
Suppose that $A\sbeq \{1,\dots, N\}$ does not contain any nontrivial three term arithmetic progression. Then
$$|A|\ll {\exp\(-c(\log N)^{1/11}\)}N.$$
\end{theorem} 
A short time later, Bloom and Sisask \cite{bs1}, by utilizing almost-periodic more efficiently in Kelley-Meka proof, were able to increase the exponent from $1/11$ to $1/9.$

However, several years earlier, even better upper bounds were already known for invariant equations equations with more  than three variables.
Shkredov and the author \cite{ss} established  such result for  invariant equations with at least six variables. 

\begin{theorem}(\cite{schoen-shkredov})
Suppose that $A\sbeq \{1,\dots, N\}$ does not contain any nontrivial solutions to the equation $x_1+\dots+x_5=5y.$ Then
$$|A|\le {\exp\(-c(\log N)^{1/7} \)}N.$$
\end{theorem}
The above result was extended in \cite{ss} to  invariant equations with at least four variables with the same upper bound.
Furthermore,  in a recent work by Ko'sciuszko \cite{ko}, it was showed, among other results, that if a subset of $\{1,\dots, N\}$ does not contain any nontrivial solutions to the equation
$x_1+\dots+x_k=ky,$ with $k\ge 2\cdot 3^{m+1}+2$, then its size is at most
$$ {\exp\(-c(\log N)^{1/(6+\g_m)} \)}N$$
where $\g_m=2^{1-m}.$ 

In this paper, we further improve upon the aforementioned results for invariant equations with at least four variables by narrowing the upper bound closer to Behrend's lower bound.

\begin{theorem}\label{main}
Suppose that $A\sbeq \{1,\dots, N\}$ does not contain any nontrivial solutions to the equation $x+y+z=3w.$ Then
$$|A|\le {\exp\(-c(\log N)^{1/5} \exp(-C\sqrt{\log\log(1/N)})\)}N.$$
\end{theorem}
We  ground our argument in a new result related to the almost-periodicity of convolutions, which is elaborated upon in Section 3 (see Theorem \ref{periodicity1}).

\bigskip
\centerline{\sc { 2. Elementary properties of Bohr sets and notation}}
\bigskip

Throughout the paper we will use Bohr sets, which are a fundamental tool
 introduced to  additive combinatorics by Bourgain \cite{b}.  Bohr sets have a rich arithmetic structure and can serve as a substitute for subspaces, which is especially useful when applying the density increment strategy.

Let $G$ be an abelian group and let us denote the dual group
of its characters by $\h G$. We define the Bohr set
with a generating set $\G\sbeq\h G$ and a radius  $\r\ge 0$
to be the set
$$B= B(\G,\r)
  =  \big\{ x\in G:  |1-\g (x)|\le\r \text { for all } \g\in\G\big\}\,.$$
	 The size of $\G$ is called the rank of $B$.
Given $\d>0$ and a
Bohr set $B=B(\G,\g)$, by $B_{\d}$ we mean the Bohr set $B(\G,\d\g).$ 
The  next lemma
is pretty standard, hence we refer the reader to \cite{tv}
for a  complete account. 

\begin{lemma}\label{bohr-size} Let $B\sbeq G$ be a Bohr set of rank $d$ and radius $\r\in [0,2]$.
 Then we have
$$
|B|\ge (\r/2\pi)^{d}N,
$$
$$|B_2|\le 6^d|B|$$
and  for every $\d\in [0,1]$
$$|B_\d|\ge (\d/2)^{3d}|B|.$$
\end{lemma}

Bohr sets do not always behave like convex bodies. The size of Bohr
sets can vary significantly even for small changes of the radius
which was the motivation behind  the following definition.
We call a Bohr set $B$ of rank $d$ and radius $\r$ {regular} if for every $\d$,
with $|\d|\le 1/100d$ we have
$$
(1-100d|\d|)|B|\le|B_{1+\d}|\le(1+100d|\d|)|B|.
$$
Bourgain \cite{b} showed that regular Bohr sets are ubiquitous.

\begin{lemma}\label{l:bohr-regular}
\label{lem:regularity} For every Bohr set $B$ there exists
$\d\in [1/2,1]$ such that  $B_\d$ is
regular.
\end{lemma}
The advantage of regularity is  that the size of sumset $B+B_\d\sbeq B_{1+\d}$ is only slightly bigger than the size of $B$, provided that $\d\le 1/100d.$

Throughout the  paper, we will use a fairly standard notation in additive combinatorics. If $T$ and $A$ are two sets, we refer  by $\m_T(A):=\frac{|A\cap T|}{|T|}$ the relative density of $A$ in $T,$ and put $\m_T:=1_T/|T|,$ where $1_T$ is the indicator function of $T$. 

For two functions $f,g:G \rightarrow \C$
we write 
$$f*g(x)=\sum_yf(y)g(x-y)$$
for the convolution of $f$ and $g$. Notice that $1_A*\m_T(x)=\frac{|A\cap (x-T)|}{|T|}$. For $p\ge 1$ define
$$\|f\|^p:=\(\sum_{x\in G} |f(x)|^p\)^{1/p},$$
and 
$$\|f\|_\infty:=\sup_{x\in G} |f(x)|.$$
For convenience, we employ the symbols $C$ and $c$ to represent positive constants that are adequately large and small, respectively. It should be noted that these constants' values may differ in various instances.  We  also use standard  Vinogradov's $\ll$ notation. By $\log$ we  always mean $\log_2.$
Moreover, when we refer to the group $G$ it is understood to be a finite abelian group; however, in  applications, we always have $G=\zn$ for a prime $N$. As usually if $A\sbeq \{1,\dots, M\}$ then we consider $A$ as a subset of $\zn$ with  $3M\le N\le 6M.$

\bigskip
\centerline{\sc { 3. A sketch of the argument}}
\bigskip

Like in all previous works concerning sets free of solutions to a convex equation, we  also employ the density increment strategy. The general line of the proof is the same as in \cite{ss}, but the details are very different. The proof method in \cite{ss} relies on the $L^\infty$-almost periodicity of three-fold convolution (see Theorem \ref{periodicity} below), which extends Sanders' observation \cite{sanders} that the Croot-Sisak \cite{c-s} $L^p$-almost periodicity lemma (for two-fold convolution) works exceptionally efficiently for the convolution $1_A*1_{A-A}.$ This allowed us in \cite{ss} to achieve a density increment at each step by a factor of $5/4$ on the Bohr set of rank increased by $\log^4(2/\a),$ where $\a$ is  the initial density of our solution-free set. 

Here we  proceed differently, applying our key ingredient, a version of $L^\infty$-almost periodicity Theorem \ref{periodicity1},  we will be able to obtain alternative result. 
This approach leads us to a result, which, in its simplest form, can be summarized as follows: either we attain a substantial increase in density within a Bohr set, whose the rank increases by $\log^{4+o(1)}(2/\alpha),$ or we achieve a modest density increment within a Bohr set, whose the rank increases  by $\log^{1+o(1)}(2/\alpha).$

Another crucial element in proving the main theorem's bound lies in the precise control of the radii of successive Bohr sets. Remarkably, Theorem \ref{periodicity1} works highly effective in this regard as well.

\bigskip
\centerline{\sc { 4. An almost-periods lemma}}
\bigskip

Here we will establish the key lemma for our approach, which is particularly useful when applying the density increment argument on a sequence of Bohr sets. The main idea is that if, after applying the following Theorem 
\ref{periodicity} proved in \cite{ss}, the shifts of set $A$ do not achieve a significant density increment on Bohr set $B,$ then we will be able to obtain the conclusion of Theorem \ref{periodicity} for a much larger Bohr set. To accomplish this, we will use a probabilistic argument, approximating the convolution $1_M*\m_B$ by a much larger set $R,$ and then applying again Theorem \ref{periodicity} to it, obtaining much more efficient estimates. The concept of applying the fact that $1_M*\m_B$ is small resembles a part the argument in \cite{schoen-roth} (specifically, refer to Lemma 14 in \cite{schoen-roth}).  Nevertheless,  this idea was employed in \cite{schoen-roth} in a completely distinct manner, making use of Fourier analysis techniques.

\begin{theorem}\label{periodicity} (\cite{ss}) Let $\e\in (0,1).$ Let $A, M, L$ be subsets of a finite abelian group $G$ and let $B\sbeq G$ be a regular Bohr set of rank $d$ and radius $\r$. Suppose $|A+S|\le K|A|$ for some subset $S\sbeq B$ with $\m_B(S)\ge \s>0,$ and assume $\eta:=|M|/|L|\le 1.$ Then there is a regular Bohr set $T$ of rank $d+d'$ and radius at least $c\r \e \eta^{1/2}/d^2d',$ where
$$d'\ll \e^{-2}\log^2(2/\e \eta)\log(2/\eta)\log(2K)+\log(1/\s)$$
such that
\begin{equation}
\|1_A*1_M*1_{L}*\m_{T}-1_A*1_M*1_{L}\|_\infty\le \e |A||M|\,.\label{l-infty}
\end{equation}
\end{theorem}

In the course of the proof of the main result of this section we will use classical Bernstein's inequality
\cite{bernstein}.

\begin{lemma}\label{bernstein} (Bernstein)
Let $ X_{1},\dots ,X_{N}$ be independent  random variables and suppose that $ |X_{k}-\E(X_k)|\leq m$ for every $1\le k\le N$. Then, for all positive $t$
$$\P\big (|\sum_{k=1}^NX_k-\sum_{k=1}^N\E(X_k)|\ge t\big )\le 2\exp\Big(-\frac{\frac12t^2}{\sum_{k=1}^N \v(X_k)+\frac13tm}\Big)\,.$$
\end{lemma}
\bigskip

\begin{theorem}\label{periodicity1} Let $\e\in (0,1).$ Let sets $A, M, L, S, B\sbeq G$   satisfy the assumptions of Theorem \ref{periodicity}. Let $T\sbeq B$ be a regular Bohr set such   that 
\begin{equation}\label{cor-1}\|1_A*1_M*1_{L}*\m_{T}-1_A*1_M*1_{L}\|_\infty\le \e |A||M|.
\end{equation}
Let  $\e_1\in (0,1)$ and suppose that for some positive $\g\le 1$ we have $1_M*\m_{T}(x)\le \g$ for every $x,$ and $|M|\ge2$. Then there is a regular Bohr set $B^1\sbeq B$ of rank $d+d_1$ and radius at least $\r \e_1 (\eta/2\g)^{1/2}/d^2d_1,$ where
$$d_1\ll \e_1^{-2}\log^2(4\g/\e_1 \eta)\log(4\g/\eta)\log(2K)+\log(1/\s)$$
such that
\begin{equation}
\|1_A*1_M*1_{L}*\m_{B_1}-1_A*1_M*1_{L}\|_\infty\le (2\e+\e_1) |A||M|+ 18|A|\sqrt{\g |M|\log |A+M+L+T|}.\label{l-infty-1}
\end{equation}
\end{theorem}

\begin{proof} First, we show that there exists a set $R\sbeq G$ such that
\begin{equation}\label{r-size}\g ^{-1}|M|-6\sqrt{\g^{-1}|M|}\le |R|\le \g ^{-1}|M|+6\sqrt{\g^{-1}|M|}
\end{equation}	
and
\begin{equation}\label{r-per}\|1_A*1_M*1_{L}*\m_{T}-\g 1_A*1_R*1_{L}\|_\infty\le  6|A|\sqrt{\g |M|\log |A+M+L+T|}.
\end{equation}	
Let  $R$ be a random subset of $G$   chosen by picking each $x\in G$ independently  with probability 
$$\P(x\in R)=\g^{-1}1_M*\m_{T}(x)\,.$$
Note that the expected size of $R$ is $\g^{-1}|M|$ and its variance does not exceed $\g^{-1}|M|$, so by Bernstein's inequality
\begin{equation}\label{r-size-1}
\P\big(\big||R|-\g^{-1}|M|\big|\le 6\sqrt{\g^{-1}|M|
}\big)\ge 3/4 \,.
\end{equation}
  Next, for a fixed $x\in W:=A+M+L+T$  the random variable 
	$$1_A*1_R*1_{L}(x)=\sum_y 1_A*1_L(x-y)1_R(y)$$ 
	has an expected value
	$$\E(1_A*1_R*1_{L}(x))=\g^{-1}1_A*1_M*1_{L}*\m_{T}(x)$$
	and its variance can be bouned from above by
	$$\v(1_A*1_R*1_{L}(x))\le \sum_y(1_A*1_L(x-y))^2\g^{-1}1_M*\m_{T}(y)\le \g^{-1}|M||A|^2.$$
	Thus, again by  Bernstein's inequality for any $x\in W$
\begin{equation*}
\P\big(\big|1_A*1_R*1_{L}(x)-\g^{-1}1_A*1_M*1_{L}*\m_{T}(x)\big|\le 6|A|\sqrt{\g^{-1} |M|\log |W|}\big)\ge 1-\frac1{4|W|} ,
\end{equation*}	
	so
\begin{equation}\label{r-convolution}
\P\big(\big\|1_A*1_R*1_{L}-\g^{-1}1_A*1_M*1_{L}*\m_{T}\big\|_\infty\le 6|A|\sqrt{\g^{-1} |M|\log |W|}\big)\ge 3/4. 
\end{equation}
Due to \eqref{r-size-1}  and \eqref{r-convolution}  there exists a set satisfying \eqref{r-size} and \eqref{r-per}.

By  Theorem \ref{periodicity} applied  with $R$ in place of $M$ and $\e_1$ in place of $\e$
there is a regular Bohr set $B^{1}\sbeq B$ of rank $d_1$ and radius at least $\r \e_1 \eta_1^{1/2}/d^2d_1,$ where
$$\eta_1=|R|/|L|\ge \frac{ \g ^{-1}|M|-6\sqrt{\g^{-1}|M|}}{|L|}\ge \eta/2\g,$$
and
$$d_1\ll \e_1^{-2}\log^2(2/\e \eta_1)\log(2/\eta_1)\log(2K)+\log(1/\s)$$
such that
\begin{equation*}
\|1_A*1_R*1_{L}*\m_{B^1}-1_A*1_R*1_{L}\|_\infty\le \e_1 |A||R|\le \e_1 \g^{-1}|A||M|+6\sqrt{\g^{-1}|M|}.\label{r-infty}
\end{equation*}
In view of \eqref{cor-1}, \eqref{r-per} and the triangle inequality we infer that
$$\|1_A*1_M*1_{L}-\g 1_A*1_R*1_{L}\|_\infty\le \e|A||M|+ 6|A|\sqrt{\g |M|\log |W|},$$
which leads to 
\begin{eqnarray*}
\|1_A*1_M*1_{L}*\m_{B^1}-\g 1_A*1_R*1_{L}*\m_{B^1}\|_\infty &\le& \|1_A*1_M*1_{L}-\g 1_A*1_R*1_{L}\|_\infty
\|\m_{B^1}\|_1 \\
&=&\|1_A*1_M*1_{L}-\g 1_A*1_R*1_{L}\|_\infty\\ 
&\le& \e|A||M|+ 6|A|\sqrt{\g |M|\log |W|}\,.\label{l-infty}
\end{eqnarray*}
By the triangle inequality we have
\begin{eqnarray*}
\|1_A*1_M*1_{L}*\m_{B^{1}}-1_A*1_M*1_{L}\|_\infty&\le& \|1_A*1_M*1_{L}-\g 1_A*1_R*1_{L}\|_\infty\\
&& + \|1_A*1_M*1_{L}*\m_{B^1}-\g 1_A*1_R*1_{L}\|_\infty \\
&\le & 
 \e |A||M|+6|A|\sqrt{\g |M|\log |W|}\\
&&+ \|1_A*1_M*1_{L}*\m_{B^1}-\g 1_A*1_R*1_{L}*\m_{B_1}\|_\infty\\
&&
+ \|\g 1_A*1_R*1_{L}-\g 1_A*1_R*1_{L}*\m_{B_1}\|_\infty\\
&\le&  (2\e+\e_1) |A||M|+18|A|\sqrt{\g |M|\log |W|},
\end{eqnarray*}
which completes the proof.$\hfill\blacksquare$
\end{proof}

\bigskip
\centerline{\sc { 4. Iterative lemmas}}
\bigskip

In this section, our goal is to establish a result (Proposition \ref{increment-main}) that will be iteratively applied in the proof of the main theorem. Similar to the approach in \cite{ss}, we divide our analysis into two cases, depending on the size of the sumset $A+A'.$ We will begin with the case where $A+A'$ is large.

\begin{lemma}\label{increment1} 
Let $B\sbeq \zn$ be a regular Bohr set of rank $d$ and radius $\r$ such that $|B|\ge (Cd/\a)^5$, and let $A\sbeq B$ has relative density $c_0\ge \m_B(A)\ge \a$ for some small constant $c_0>0$. Let $B'=B_\d$,  where $\d=1/Cd$, be a regular Bohr set such that $B_{1+3\d}\le 1.01|B|$ and assume that $A'=A\cap B$ satisfies $\m_{B'}(A')\ge \a$ and $|A+A'|\ge
|A|/2\a.$ Let $2\le  h\le \log(1/\a)$ be a real number.  If  $A$ does not contain any nontrivial solutions  to $x+y+z=3w$, then there is a positive integer $k\le  \lceil\log\log(1/\a^8)/\log h \rceil$ and a Bohr set $T\sbeq B$ of rank at most $d+d'$ and radius  at least 
$$c\r \a^{1/2h^{k-1}}/d^4d'\log (1/\a)^{1/\log h},$$
 where 
$$d'\ll \(\log \log (1/\a)+h^{-(k-1)}\log (1/\a)\)^3\(\log (1/\a)\)^{1+2/\log h}$$
 such that
$\|1_A*\m_T\|_\infty\ge \frac{10}7\a^{1-1/h^{k}}.$
\end{lemma}
\begin{proof}  We apply Theorem \ref{periodicity} and Theorem \ref{periodicity1} with $A:=-3\cdot A',$  
$S:=3\cdot B'_\n$, where $\n\le 1/Cd$, $M:=A$ and $B_{1+3\d}\setminus (A+A')$ in place of $L$, and 
$$\e=\frac1{160}(\log(1/\a^8))^{-1/\log h}.$$
 Observe that by regularity of $B'$ we have
\begin{equation*}
|3\cdot A'+S|\le |B'_{1+\n}|\le \frac{2}{\a}|3\cdot A'|.
\end{equation*}
By Theorem \ref{periodicity} there is a Bohr set $B^1\sbeq S$ of rank at most $d+d_1$ and radius at least $c\r\e\a^{1/2}/d^4d_1$, where
$d_1\ll \e^{-2}\log^4(1/\a)\ll \(\log (1/\a)\)^{4+2/\log h}$ such that
$$\|1_{-3\cdot A'}*1_A*1_{L}*\m_{B^1}-1_{-3\cdot A'}*1_A*1_{L}\|_\infty\le \e |A'||A|.$$
If  $\|1_A*\m_{B^1}\|_\infty\ge \frac{10}7\a^{1-1/h}$ the proof is concluded; otherwise, we will  apply Theorem \ref{periodicity1} with $\e_1=\e$ and $\g=\frac{10}{7}\a^{1-1/h}$. 
Then there is a Bohr set $B^2\sbeq S$
 of rank $d+d_2$  and radius at least 
$$c \r \e\a^{1/2h}/d^4d_2,$$
 where
\begin{eqnarray*}
d_2&\ll& \e^{-2}\log^2(14/5\e\a^{1/h})\log(14/5\a^{1/h})\log (2/\a)\\
&\ll& C\e^{-2}\log^2(3/\e\a^{1/h})\log(3/\a^{1/h})\log (2/\a)\\
&\ll& \(\log\log (1/\a)+h^{-1}\log (1/\a)\)^3\(\log (1/\a)\)^{1+2/\log h}
\end{eqnarray*}
 such that
\begin{eqnarray*}
\|1_{-3\cdot A'}*1_A*1_{L}*\m_{B^2}-1_{-3\cdot A'}*1_A*1_{L}\|_\infty &\le& 3\e |A'||A|\\
&&+18|A'|\sqrt{\frac{10}{7}\a^{1-1/h} |A|\log |A+L-3\cdot A'+B^1|}\\ 
&\le& 3\e |A'||A|+18|A'|\sqrt{ |A|\log ( 36^{d}|B|)},
\end{eqnarray*}
the last inequality results from the assumption that $\a\le c_0$, and Lemma \ref{bohr-size}  
$$|A+L-3\cdot A'+B^1|\le |B_{3+3\d+3\n}|\le |B_4|\le 36^{d}|B|.$$
We will apply the above procedure iteratively, and suppose that after $k$ steps we obtain a Bohr set $B^k\sbeq S$ of rank $d+d_k$
and radius 
$$c\r \e\a^{1/2h^{k-1}}/d^{4}d_k,$$ where 
$$d_k\ll  \(\log \log (1/\a)+h^{-(k-1)}\log (1/\a)\)^3\(\log (1/\a)\)^{1+2/\log h}$$
 such that
\begin{eqnarray*}
\|1_{-3\cdot A'}*1_A*1_{L}*\m_{B^k}-1_{-3\cdot A'}*1_A*1_{L}\|_\infty&\le& (2^{k}-1)\e |A'||A|\\
&&+
(2^{k-1}-1)18|A'|\sqrt{ |A|\log |A+L-3\cdot A'+B^{k-1}|}\\
&\le& (2^{k}-1)\e |A'||A|+
(2^{k-1}-1)18|A'|\sqrt{ |A|\log ( 36^{d}|B|)},
\end{eqnarray*}
due to the inequality $|A+L-3\cdot A'+B^{k-1}|\le |B_4|.$
If $\|1_A*\m_{B^k}\|_\infty\ge \frac{10}7\a^{1-1/h^{k}}$ the proof is concluded; otherwise, we will  apply Theorem \ref{periodicity1} with $\e_1=\e$ and $\g=\frac{10}7\a^{1-1/h^{k}}.$ Thus, there is a Bohr set $B^{k+1}\sbeq S$ of rank $d+d_{k+1}$ and radius 
$$c\r \e\a^{1/2h^{-k}}/d^{4}d_{k+1},$$ 
where 
$$d_{k+1}\ll C\(\log \log (1/\a)+h^{-k}\log (1/\a)\)^3\(\log (1/\a)\)^{1+2/\log h}$$
 such that
$$\|1_{-3\cdot A'}*1_A*1_{L}*\m_{B^{k+1}}-1_{-3\cdot A'}*1_A*1_{L}\|_\infty\le (2^{k+1}-1)\e |A'||A|+(2^k-1)18|A'|\sqrt{ |A|\log ( 36^{d}|B|)}.$$
If for every $k< l:=\lceil\log\log(1/\a^8)/\log h \rceil $ we have $\|1_A*\m_{B^k}\|_\infty< \frac{10}{7}\a^{1-1/h^{k}}$  then again by  Theorem \ref{periodicity1}
there is a Bohr set $B^{l}\sbeq S$ of rank $d+d_l$ 
and radius 
$$c\r\e \a^{1/2h^{l-1}}/d^4d_{l},$$
 where 
$$d_l\ll\(\log \log (1/\a)+h^{-(l-1)}\log (1/\a)\)^3\log(2/\a)\(\log (2/\a)\)^{1+2/\log h}$$ 
such that
\begin{eqnarray*}
\|1_{-3\cdot A'}*1_A*1_{L}*\m_{B^{l}}-1_{-3\cdot A'}*1_A*1_{L}\|_\infty &\le& 
(2^{l}-1)\e |A'||A|+(2^{l-1}-1)18|A'|\sqrt{ |A|\log ( 36^{d}|B|)}\\
&\le&\frac1{80} |A'||A|+(2^l-1)18|A'|\sqrt{ |A|\log ( 36^{d}|B|)}\\
&\le& \frac1{80} |A'||A|+20\log(1/\a)|A'|\sqrt{ |A|\log ( 36^{d}|B|)}.
\end{eqnarray*}
Next, we show that the last term of the right hand side of the above inequality does not exceed 
$\frac1{80} |A'||A|$. Note that the inequalities
$$ 400(\log(1/\a))^25d\le 2000\a^{-2}d\le \frac1{160^2}\a|B|\le \frac1{160^2}|A|$$
and
$$400(\log(1/\a))^2\log |B|\le \frac1{160^2}|A|$$
hold provided that $|B|\ge (Cd/\a)^5$ and $C$ is large enough,
hence 
$$20\log(1/\a)|A'||A|^{1/2}\log ^{1/2}( 64^{d+1}|B|)\le \frac1{80} |A'||A|$$
and therefore we have
$$\|1_{-3\cdot A'}*1_A*1_{L}*\m_{B^{l}}-1_{-3\cdot A'}*1_A*1_{L}\|_\infty\le \frac1{40} |A'||A|.$$
The last step is treated differently as we now utilize  the fact that $A$ is a solution free set. Since there are only trivial solutions in $A$ to $x+y+z=3w$ it follows that
$$1_{-3\cdot A'}*1_A*1_{L}(0)=1_{-3\cdot A'}*1_A*1_{B_{1+3\d}}(0)-1_{-3\cdot A'}*1_A*1_{A+A'}(0)=|A'||A|-|A'|.$$
Thus
\begin{equation*}
1_{-3\cdot A'}*1_A*1_{L}*\m_{B^{l}}(0)\ge |A'||A|-|A'|-\frac1{40} |A'||A|\ge \frac{19}{20}|A'||A|,
\end{equation*}
as $A$ is large enough, hence
$$\frac{19}{20}|A'||A|\le |A'||L|\|1_A*\m_{B^k}\|_{\infty}\le \frac{0.501}{\a}|A'||A|\|1_A*\m_{B^l}\|_{\infty},$$
so
$$\|1_A*\m_{B^l}\|_{\infty}\ge \frac95\a.$$
To complete the proof it is enough to observe that
$$\a^{-1/h^{l}}\le \a^{-\frac1{8\log(1/\a)}}=2^{1/8}<1.1$$ 
hence
$$\|1_A*\m_{B^l}\|_{\infty}\ge\frac95\a\ge \frac{10}7\a^{1-1/h^{l}}$$
and the assertion follows. $\hfill\blacksquare$
\end{proof}

\bigskip

The density increment in the case $|A+A'|\le |A|/2\a$ will be proved by a similar reasoning, so we  present the argument in a somewhat condensed form. Here, we do not  even need the assumption that $A$ does not contain any nontrivial solutions  to the equation $x+y+z=3w.$

\begin{lemma}\label{increment2}
Let $B\sbeq \zn$ be a regular Bohr set of rank $d$ and radius $\r$ such that $|B|\ge (Cd/\a)^5$, and let $A\sbeq B$ has relative density $c_0\ge \m_B(A')\ge \a$, for some small constant $c_0>0$, and suppose that 
$|A+A'|\le |A|/2\a.$  Let $2\le  h\le \log(1/\a)$ be a real number. Then there is a positive integer $k\le  \lceil\log\log(1/\a^8)/\log h \rceil $ and a Bohr set $T\sbeq B$ of rank at most $d+d'$ and radius  at least 
$$c\r \a^{1/2h^{k-1}}/d^4d'\log (1/\a)^{1/\log h},$$
 where 
$$d'\ll \(\log \log (1/\a)+h^{-(k-1)}\log (1/\a)\)^3\(\log (2/\a)\)^{1+2/\log h}$$
 such that $\|1_A*\m_{B^k}\|_\infty\ge \frac{10}7\a^{1-1/h^{k}}.$
\end{lemma}
\begin{proof} Set $S=B_\n$, where $\n:=1/Cd$, then by regularity of $B$ we have
$$|A'+S|\le |B_{1+\n}|\le \frac2{\a}|A'|.$$
  We apply Theorem  \ref{periodicity} and Theorem \ref{periodicity1} with $A:= A',$ the set $S$,  $M:=A$, $L:=-A-A'$, and 
	$$\e=\frac1{40}(\log(1/\a^8))^{-1/\log h}$$
	to get a Bohr set  $B^1\sbeq S$ of rank at most $d+d_1$ and radius at lest $c\r\e\a^{1/2}/d^3d_1$, where
$d_1\ll \e^{-2}\log^4(1/\a)\ll \(\log (1/\a)\)^{4+2/\log h}$ such that
$$\|1_{ A'}*1_A*1_{L}*\m_{B^1}-1_{ A'}*1_A*1_{L}\|_\infty\le \e|A'||A|.$$
If $\|1_A*\m_{B^1}\|_\infty\ge \frac{10}7\a^{1-1/h}$  the proof is concluded; otherwise, we will  apply Theorem \ref{periodicity1} with $\e_1=\e$ and $\g=\frac{10}{7}\a^{1-1/h}$ to get
 a Bohr set $B^2$
 of rank $d+d_2$, where  and radius at least 
$$c\r \e \a^{1/2h}/d^3d_2,$$ 
where
$$d_2\ll \(\log \log (1/\a)+h^{-1}\log (1/\a)\)^3\(\log (2/\a)\)^{1+2/\log h}$$ 
such that
$$\|1_{A'}*1_A*1_{L}*\m_{B^2}-1_{A'}*1_A*1_{L}\|_\infty\le 3\e |A'||A|+18|A'|\sqrt{ |A|\log ( 36^{d}|B|)}.$$
After $k$ iterations we obtain a Bohr set $B^k$ of rank $d+d_k$ and radius 
$$c\r \e\a^{1/2h^{k-1}}/d^3d_k,$$
where 
$$d_k\ll \(\log \log (1/\a)+h^{-(k-1)}\log (1/\a)\)^3\(\log (2/\a)\)^{1+2/\log h},$$ 
such that
$$\|1_{ A'}*1_A*1_{L}*\m_{B^k}-1_{A'}*1_A*1_{L}\|_\infty\le (2^k-1)\e |A'||A|+(2^{k-1}-1)18|A'|\sqrt{ |A|\log ( 36^{d}|B|)}.$$
We will repeat this argument  unless we achieve required density increment, but no more than  $l:= \lceil\log\log(1/\a^8)/\log h \rceil$ times. 
If we do not obtain  required density increment for any $k<l$, then there is a 
there is a Bohr set $B^{l}$ of rank $d+d_l$ and radius 
$$c\r \e\a^{1/2h^{l-1}}/d^3d_{l},$$
 where 
$$d_l\ll \(\log \log (1/\a)+h^{-(l-1)}\log (1/\a)\)^3\(\log (1/\a)\)^{1+2/\log h}$$ 
such that
$$
\|1_{ A'}*1_A*1_{L}*\m_{B^{l}}-1_{A'}*1_A*1_{L}\|_\infty\le \frac1{10}|A'||A|.
$$
Clearly, $1_{A'}*1_A*1_{L}(0)=|A||A'|$, so 
$$1_{ A'}*1_A*1_{L}*\m_{B^{l}}(0)\ge  |A||A'|-\frac1{20} |A'||A|\ge \frac{9}{10}|A||A'|,$$
hence
$$\|1_A*\m_{B^l}\|_{\infty}\ge \frac95\a\ge  \frac{10}7\a^{1-1/h^l}.$$
This concludes the proof.$\hfill\blacksquare$
\end{proof}

\bigskip

The next lemma is quite standard; however, we cannot utilize its analogous version proven in \cite{ss} 
(Lemma 6.4). This is because in \cite{ss}, it is proven for much smaller $\d\le \a/Cd$, which is insufficient for our approach and would not allow us to increase the exponent to $1/5$ in the main result.
It turned out, however, that only a minor modification of the proof allows us to show the same thesis for a significantly larger  $\d\le 1/Cd$.

\begin{lemma}\label{density} Let $B$ be a regular Bohr set of rank $d$, let $A\sbeq B$ has relative density $\m_B(A)=\a.$ Let $C>0$ be a constant such that $|B|\ge (1-\frac1{20}) |B_{1+\d}|,$ $\d=1/Cd$ 
and let $B', B''\sbeq B_\d.$ Then either:
\begin{enumerate}
\item there is an $x$ such that $1_A*\m_{B'}(x)\ge \frac8{10}\a$ and $1_A*\m_{B''}(x)\ge \frac8{10}\a$; or
\item $\|1_A*\m_{B'}\|_\infty\ge 1.1\a$, or $\|1_A*\m_{B''}\|_\infty\ge 1.1\a.$
\end{enumerate}
\end{lemma}
\begin{proof} 
Since $A+B', A+B''\sbeq B_{1+\d}$ it follows that
$$\sum_{x\in B_{1+\d}}1_A*\m_{B'}(x)=\sum_{x\in B_{1+\d}}1_A*\m_{B''}(x)=|A|,$$
hence
$$\sum_{x\in B_{1+\d}}\(1_A*\m_{B'}(x)+1_A*\m_{B''}(x)\)=2|A|.$$
Thus, for some $x$ we have 
$$1_A*\m_{B'}(x)+1_A*\m_{B''}(x)\ge 2|A|/|B_{1+\d}|=2\a|B|/|B_{1+\d}|\ge \(2-\frac1{10}\)\a.$$
If $\|1_A*\m_{B'}\|_\infty, \|1_A*\m_{B''}\|_\infty <1.1\a$ then 
$$1_A*\m_{B'}(x),1_A*\m_{B''}(x)\ge \(2-\frac1{10}\)\a-1.1\a= \frac8{10}\a,$$
and the proof is completed.
$\hfill\blacksquare$
\end{proof}
\bigskip

\no Now, we are in position to prove the main result of this section.

\begin{proposition}\label{increment-main}
Let $B\sbeq \zn$ be a regular Bohr set of rank $d$ and radius $\r$, and let $A\sbeq B$ has relative density $c_0\ge \m_B(A)\ge \a.$  Assume that $|B|\ge (Cd/\a)^{5d}.$ Let $2\le  h\le \log(2/\a)$ be a real number. If  $A$ does not contain any nontrivial solutions  to $x+y+z=3w$ then there  exist a positive integer $k\le \lceil\log\log(1/\a^9) /\log h\rceil$ and a Bohr set $T\sbeq B$ of rank at most $d+d'$ and radius  at least 
$$c\r \a^{1/2h^{k-1}}/d^5d'\(\log (1/\a)\)^{1/\log h},$$
 where 
$$d'\ll \(\log \log (1/\a)+h^{-(k-1)}\log (1/\a)\)^3\(\log (1/\a)\)^{1+2/\log h}$$
such that $\|1_A*\m_T\|_\infty\ge \a^{1-1/h^{k}}.$
\end{proposition}
\begin{proof} 
We chose constants  $c',c''\in [1/2, 1]$ such that the Bohr sets $B'=B_{\d'}$  and $B''=B'_{\d''},$ where $\d'=c'/100d$ and $\d''=c''/100d$, are regular and $|B_{1+3\d}|\le 1.01 |B|.$ 
If the second conclusion of Lemma \ref{density} holds then we have 
$$\|1_A*\m_{B'}\|_\infty\ge 1.1\a \text{ or } \|1_A*\m_{B''}\|_\infty\ge 1.1\a.$$
Clearly, the Bohr sets $B'$ and $B''$ have rank $d$ and radius at least $c\r/d^2.$ In this case we 
can take $k= \lceil \log\log(1/\a^8)/\log h\rceil$  to get the required inequalities. Indeed, for such choice of $k$ we have 
$$\a^{1-1/h^{k}}\le 2^{1/9}\a<1.1\a.$$
If the first conclusion of Lemma \ref{density} holds then for some $x\in B$ we have
$$a_1:=1_A*\m_{B'}(x)\ge \frac8{10}\a \text{~ and ~} \a_2:=1_A*\m_{B''}(x)\ge \frac8{10}\a,$$
and  by the above, we can assume that
$\a_1,\a_1\le 1.1\a.$ 
In order to apply Lemma \ref{increment1} and Lemma \ref{increment2}, it is necessary to ensure that both $\alpha_1$ and $\alpha_2$ do not exceed $c_0.$
If this is not the case, we find a subset $A_1\sbeq A$ such that
$$c_0\ge 1_{A_1}*\m_{B'}(x)\ge \frac7{10}\a \text{~ and ~} c_0\ge 1_{A_1}*\m_{B''}(x)\ge \frac7{10}\a.$$ 
  Since we have $|B'|\ge 2$ by assumption, the Cauchy-Davenport theorem yields
\begin{equation}\label{B'}
|B'|\ge |\lfloor 100d/c''\rfloor B''|\ge \lfloor 100d/c''\rfloor|B''|-\lfloor 100d/c''\rfloor+1\ge \frac {\lfloor 100/c''\rfloor}2|B''|\ge 25|B''|.
\end{equation}
 Let $U'\sbeq (A-x)\cap (B'\setminus B'')$ and $U''\sbeq (A-x)\cap B''$ be arbitrary subsets such that 
$|U''|=\max\(0, (\a_2-c_0)|B''|\)$ and $|U'|=\max\(0, (\a_1-c_0)\)|B'|-|U''|.$ We show that $A_1$ can be taken as $A\setminus (U'\cup U'').$ Indeed, we have
$$c_0\ge 1_{A_1}*\m_{B''}(x)=\a_2-\max\(0, (\a_2-c_0)\)\ge \frac8{10}\a$$
and due to \eqref{B'}
\begin{eqnarray*}
c_0\ge 1_{A_1}*\m_{B'}(x)&=&\a_1-\max\(0, (\a_1-c_0)\)-\frac1{25}\(\a_2-\max\(0, (\a_2-c_0)\)\)\\
&\ge& \frac8{10}\a-\frac1{25}1.1\a\ge \frac7{10}\a.
\end{eqnarray*}
We put $A':=(A_1-x)\cap B', A'':=(A_1-x)\cap B''$ and observe that $A'$ is solution free to the equation $x+y+z=3w.$

If $|A'+A''|\ge |A'|/2\l,$  where $\l:=\frac7{10}\a$ then by Lemma \ref{increment1} applied with $$2\le  h\le {\log(1/\l)}\le  {\log(2/\a)}$$ there is a Bohr set $T$ of rank $d+d'$ and radius at least 
$$c(\r/C'd)\l^{1/2h^{k-1}}/d^4d'\(\log (1/\l)\)^{1/\log h}\gg \r\a^{1/2h^{k-1}}/d^5d'\(\log (1/\a)\)^{1/\log h},$$ 
 with 
\begin{eqnarray*}
d' &\ll& \(\log \log (1/\l)+h^{-(k-1)}\log (1/\l)\)^3\(\log (1/\l)\)^{1+2/\log h}\\
&\ll& \(\log \log (1/\a)+h^{-(k-1)}\log (1/\a)\)^3\(\log (1/\a)\)^{1+2/\log h}
\end{eqnarray*}
 for some $k\le \lceil \log\log(1/\l^8)/\log h\rceil \le \lceil \log\log(1/\a^9)/\log h\rceil$ such that 
$$\|1_A*\m_T\|_\infty\ge \|1_{A'}*\m_T\|_\infty\ge \frac{10}7\l^{1-1/h^{k}}\ge \a^{1-1/h^{k}}.$$

If $|A'+A''|\le |A'|/2\l,$  then by Lemma \ref{increment2} there is a Bohr set $T$ of rank $d+d'$ and radius at least 
$$c(\r/C'C''d^2)\a^{1/2h^{k-1}}/d^3d'\(\log (1/\a)\)^{1/\log h}\gg \r\a^{1/2h^{k-1}}/d^5d'\(\log (1/\a)\)^{1/\log h},$$  with 
$$d'\le C\(\log \log (1/\a)+h^{-(k-1)}\log (1/\a)\)^3\(\log (2/\a)\)^{1+2/\log h}$$
 for some $k\le  \lceil \log\log(1/\l^8)/\log h\rceil \le \lceil \log\log(1/\a^9)/\log h\rceil$ such that 
$$\|1_A*\m_T\|_\infty\ge \|1_{A'}*\m_T\|_\infty\ge \frac{10}7\l^{1-1/h^{k}}\ge \a^{1-1/h^{k}}.$$
  $\hfill\blacksquare$

\bigskip

\no{\it Remark} Notice  that the density increment obtained in Proposition  \ref{increment-main}  is always at least 
\begin{equation}\label{h-increment}
\a^{-1/h^{l}}\ge 2^{1/9h}.
\end{equation}
\bigskip

Let us  also refer to a similar lemma from [3]. Despite implying weaker result than Proposition 
\ref{increment-main}, we will also apply it in the case of $\a>c_0.$

\begin{proposition}\label{increment-old}
Let $B\sbeq \zn$ be a regular Bohr set of rank $d$ and radius $\r$, and let $A\sbeq B$ has relative density $\m_B(A)\ge \a.$ Assume that $|B|\ge (Cd/\a)^{3d}.$ If  $A$ does not contain any nontrivial solutions  to $x+y+z=3w$ then there exist  a Bohr set $T\sbeq B$ of rank at most $d+d'$ and radius  at least $\r \a^{3/2}/d^5d'$, where 
$d'\ll \log ^4(2/\a)$
such that $\|1_A*\m_T\|_\infty\ge \frac54\a.$
\end{proposition}

\end{proof}
\bigskip
\centerline{\sc { 2. Proof of Theorem \ref{main}}}
\bigskip

Now, we  iteratively apply the Proposition \ref{increment-main} and Proposition \ref{increment-old} to prove Theorem \ref{main}.
We start with $A_0=A, B^0=\zn,$ so $d_0=1, \r_0 =2$ and $\a_0=\a.$ We will apply the following iteration
scheme as long as it is possible.

If $|B^i|\ge (Cd_i/\a_i)^{5d_i}$ and $\a_i\le c_0,$  then we apply Proposition \ref{increment-main} with  
\begin{equation}\label{hi}
h_i:=\exp(\sqrt{\log\log (1/\a_i)})
\end{equation}
 to get a Bohr set $B^{i+1}\sbeq B^i$ of rank $d_{i+1},$
radius $\r_{i+1}$ and a positive integer 
$$k_{i}\le \lceil \log\log(2/\a_i)/\log h_i \rceil$$
 such that
\begin{equation}\label{d}
d_{i+1}\le d_i+C\(\log \log (1/\a_i)+h_i^{-(k_{i}-1)}\log (1/\a_i)\)^3\(\log (2/\a_i)\)^{1+2/\log h_i},
\end{equation}
\begin{equation}\label{r}
\r_{r+1}\ge c\r_i\a_i^{1/2h_i^{(k_{i}-1)}}/d_{i}^5d_{i+1}(\log 1/\a_i)^{1/\log h_i},
\end{equation}
and a set $A_{i+1}=(A_i-x_i)\cap B^{i+1}$ of relative density on $B^{i+1}$
\begin{equation}\label{ki}
\a_{i+1}\ge \a_i^{1-1/h_i^{k_i}}.
\end{equation}
Let us here recall that our equation is translation invariant, so $(A_i-x)\cap B^{i+1}$ for any $x$ is still free of solutions to our equation.

If $|B^i|\ge (Cd_i/\a_i)^{5d_i}$ but $\a_i> c_0,$ then we apply Proposition \ref{increment-old} to get a Bohr set $B^{i+1}\sbeq B^i$ of rank $d_{i+1}$
 and radius $\r_{i+1}$  such that
\begin{equation}\label{d1}
d_{i+1}\le d_i+C\log^4 (2/\a_i) ,
\end{equation}
\begin{equation}\label{r1}
\r_{r+1}\ge c\r_i\a_i^{3/2}/d_{i}^5d_{i+1},
\end{equation}
and  a set $A_{i+1}=(A_i-x_i)\cap B^{i+1}$ of relative density on $B^{i+1}$
\begin{equation}\label{ki1}
\a_{i+1}\ge \frac54\a_i.
\end{equation}
Note that because $\a_i>c_0$, we can apply Proposition \ref{increment-old} only a constant number of times.

Since the density is naturally bounded from above by $1$ and the growth of $\a_i$, by \eqref{h-increment} and \eqref{hi}, in each step is at least by the factor
$$2^{\frac1{9\exp(\sqrt{\log\log (1/\a)})}},
$$ 
so after  
$$s\ll\log (1/\a)\exp(\sqrt{\log\log (1/\a)})$$
 iterations we will be not able to continue this process. This implies that the condition $|B^s|\ge (Cd/\a_s)^{5d_s}$ must be violated. Let $t$ be the number of steps where we applied Proposition \ref{increment-main}, clearly  $s-t=O(1)$. Thus, by Lemma \ref{bohr-size}
we have
$$\r_s^{d_s} N< (Cd_s/\a_s)^{5d_s},$$
so 
\begin{equation}\label{final} N<  (Cd_s/\r_s\a_s)^{5d_s}
\end{equation}
Observe that by \eqref{ki}
$$\a^{\prod_{i=0}^{t-1}(1-1/h_i^{k_i})}\le c_0,$$
hence
$$\prod_{i=0}^{t-1}(1-1/h_i^{k_i})\gg \log^{-1}(1/\a).$$
Together with the inequality $1-1/x<e^{-x}$ this yields
\begin{equation}\label{h}
\sum_{i=0}^{t} 1/h_i^{k_i}\ll \log\log(1/\a).
\end{equation}
Since 
$$1/h_i\ge \exp(-\sqrt{\log\log(1/\a)})$$
 for every $i$, it follows from \eqref{h} that
\begin{equation*}
\sum_{i=0}^{t} 1/h_i^{k_i}\ge  \exp(-2\sqrt{\log\log(1/\a)})\sum_{i=0}^{t} 1/h_i^{k_i-1},
\end{equation*}
so
\begin{equation}\label{h-1}
\sum_{i=0}^{t} 1/h_i^{k_i-1}\ll \exp(2\sqrt{\log\log(1/\a)}).
\end{equation}
Therefore, by  \eqref{d}, \eqref{d1} and \eqref{h-1} we infer that
\begin{eqnarray*}
d_s&\ll& \sum_{i=0}^{t}\(\log \log (1/\a_i)+h_i^{-(k_{i}-1)}\log (1/\a_i)\)^3(\log^4(1/\a_i))^{4+2/\log h_i}+(s-t)\log^4(2/\a)\\
&\ll& t\(\log\log(1/\a)\)^3+\sum_{i=0}^{t}h_i^{-3(k_i-1)}(\log^4(1/\a_i))^{4+2/\log h_i}+(s-t)\log^4(2/\a)\\  
&\ll& \log^4(2/\a)\exp(C\sqrt{\log\log(1/\a)}).
\end{eqnarray*}
In view of  \eqref{r}, \eqref{r1} and \eqref{h-1} we have
\begin{eqnarray*}
\r_s&\ge& c^s\a^{\sum_{i=0}^t1/2h_i^{k_i-1}+\frac32(s-t)}/(\prod d_i^6)(\log 1/\a_i)^{s/\log h_i}\\
&\ge& c^s\a^{ C\exp(2\sqrt{\log\log(1/\a)})} / (\log(1/\a ))^{C\log(1/\a)}\\
&\ge& c\a^{C' \exp(2\sqrt{\log\log(1/\a)})}.
\end{eqnarray*}
Inserting the above estimates to \eqref{final} we get
$$N\le \exp(C\log^5(1/\a) \exp(C\sqrt{\log\log(1/\a)}),$$
and the assertion follows.

\bigskip
\centerline{\sc { 3. Concluding remarks}}
\bigskip

 It is conjectured in \cite{green-tao} that if  $\m_G(A)=\a$ then $2A-2A$ contains a Bohr set of rank $C\log(1/\a)$ and radius $\a^C$. Theorem \ref{periodicity} is closely related to  Bogolyubov-type results, as seen in \cite{sanders}. The currently best estimate of  rank of the Bohr set in Bogolyubov's lemma is due to Sanders \cite{sanders} and is bounded from above by $C\log^4(2/\a)$. However, any improvement in Theorem \ref{periodicity}  leads to a better estimate in Bogolyubov's lemma. 
This provides some support for the conjecture that in Theorem \ref{periodicity}, we can replace roughly $ C\log^4(2/\a)$ with $C\log(2/\a)$. Assuming this conjecture is true and using our method we would obtain the upper bound  
$$|A|\ll \exp\(-c (\log N)^{1/2-o(1)}\)N$$
for any set $A\sbeq \{1,\dots, N\}$ without nontrivial solutions to $x+y+z=3w,$ which essentially 
matches Behrend's lower bound. 

Since Theorem \ref{periodicity} is also a crucial component of the proof of the Kelley-Meka theorem, a natural question arises whether the argument presented in this work, based on Theorem \ref{periodicity1},  can be applied in the proof of the Kelley-Meka theorem. However, Theorem  \ref{periodicity} is not employed directly on the set $A$ in the proof of the Kelley-Meka theorem. Instead, it is applied to sets of the form $X = (A+s_1)\cap \dots \cap (A+s_p),$ which have significantly lower density. The application of Theorem \ref{periodicity1} would only lead to a relatively minor increase
of the density of the set $X,$ which seems to be of marginal relevance.

\end{document}